\newcommand{\RP}{{\mathbb RP}}
\newcommand{\Z}{{\mathbb Z}}
\newcommand{\R}{{\mathbb R}}

\newcommand{\sm}{\setminus}
\newcommand{\F}{{\mathcal F}}

\documentclass[12pt]{amsart}

\begin{document}

\author{V.A.~Vassiliev}
\thanks{
Supported by RFBR (project 95-01-00846a), INTAS (project 4373) and Netherlands
Organization for Scientific Research (NWO), project 47.03.005.}

\title{Homology of spaces of homogeneous polynomials
in $\R^2$ without multiple zeros}

\date{Revised version published in 1998}

\maketitle

\section{Introduction}

For any natural $d \ge k \ge 2$ we calculate the cohomology groups of the space
of homogeneous polynomials $\R^2 \to \R$ of degree $d$, which do not vanish
with multiplicity $\ge k$ on real lines. For $k=2$ this problem provides the
simplest example of the situation, when the ``finite-order'' invariants of
nonsingular objects are not a complete system of invariants.

The ``affine'' version of this problem (the calculation of the homology group
of the space of polynomials $\R^1 \to \R^1$ with leading term $x^d$ and without
$\ge k$-fold roots) was solved by V.~I.~Arnold in \cite{A89}, see also
\cite{V89}. As in these works, our present calculation is based on the study of
the discriminant set, i.e. of the set of polynomials with forbidden multiple
zeros.

The problem solved below turns out to be more complicated, because an essential
contribution to the homology group comes from the neighborhood of the
``infinitely degenerate'' polynomial equal identically to 0. By this reason,
the method of simplicial resolutions of the discriminant set, solving
immediately the ``affine'' problem, is replaced by its continuous analog:
conical resolution, used previously in \cite{V91}.

Here we have the simplest situation, when the invariants of ``finite order'' of
the space of nonsingular objects do not constitute a complete system of
invariants. Indeed, our spaces of nonsingular polynomials can be considered as
finite-dimensional approximations of the space $\F \sm \Sigma_k$ of smooth
functions $S^1 \to \R^1$ without $k$-fold zeros (if $d$ is odd, then with
values in a nontrivial line bundle); here $S^1$ is realized as a half of the
unit circle in $\R^2.$ By analogy with \cite{V90}, the cohomology classes of
``finite order'' of the space $\F \sm \Sigma_k$ are exactly those, which are
obtained by a natural stabilization of such cohomology groups for approximating
spaces. It turns out, that for $k=2$ and even $d$ all such 0-dimensional
cohomology classes (i.e. the invariants of degree $d$ polynomial functions $S^1
\to \R^1$ without multiple zeros) are polynomials of the number of zeros of
these functions, in particular they cannot separate everywhere positive
functions from the everywhere negative ones. For any finite even $d$ such
polynomial functions can be distinguished by a certain 0-dimensional cohomology
class, arising from the construction of the conical resolution, but this class
is not stable.

Note that a similar example, proving that the system of all finite-order
invariants of knots in $\R^3$ (see \cite{V90}) is not complete, is not
constructed yet.

\section{Main result}

{\sc Notation.} Denote by $HP_d$ the space of homogeneous polynomials $\R^2 \to
\R^1$ of degree $d$, and by  $\Sigma_k \subset HP_d$ the set of polynomials
taking zero value with multiplicity $\ge k,$ $k \ge 2,$ on some line in $\R^2$.
For any topological space $X$, $B(X,j)$ denotes its $j$-th configuration space
(i.e. the space of $j$-point subsets in $X$, supplied with a natural topology).
\medskip

{\sc Main Theorem.} {\it 1. If $k$ is even, then the group $\tilde H^*(HP_d \sm
\Sigma_k)$ is free Abelian of rank $2[d/k]+1$, and its free generators have
dimensions $k-2, 2(k-2), \ldots, [d/k](k-2),$ $k-1, 2(k-2)+1, \ldots
,[d/k](k-2)+1$ and $d-2[d/k].$

2. If $k$ is odd and $d$ is not a multiple of $k$, then the group $\tilde
H^*(HP_d \sm \Sigma_k)$ is the direct product of following groups:

(a) for any $p= 1, 2, \ldots, [d/k]$ such that $d-p\cdot k$ is odd, $\Z$ in
dimension $p(k-2)$ and $\Z$ in dimension $p(k-2)+1$;

(b) for any $p= 1, 2, \ldots, [d/k]$ such that $d-p \cdot k$ is even, $\Z_2$ in
dimension $p(k-2)+1$;

(c) $\Z$ in dimension $d-2[d/k]$.

3. If $k$ is odd and $d$ is a multiple of $k$, the answer is almost the same as
in the case 2, only the summand $\Z_2$ in dimension $d-2(d/k)+1$ vanishes.}
\medskip

{\bf Main example.} Let $d$ be even and $k=2$. Then the space $HP_d \sm
\Sigma_k$ consists of $d/2+2$ connected components, two of which (corresponding
to everywhere positive and everywhere negative functions) are contractible, and
all the other are homotopy equivalent to a circle. This homotopy equivalence is
a composition of two: the first maps any polynomial to the collection of its
zero lines (i.e. to an element of the configuration space $B(\RP^1,2p)$ with
appropriate $p \in [1,d/2]$), and the second is the arrow in the following
well-known statement.
\medskip

{\sc Lemma 1.} {\it For any $j,$ there is a locally trivial fibre bundle
$B(S^1,j) \to S^1$, whose fiber is homeomorphic to an open $(j-1)$-dimensional
disc. This fibre bundle is trivial if $j$ is odd and is non-orientable if $j$
is even.} \quad $\Box$ \medskip

{\bf Caratheodory theorem.} In the proof of the main theorem we use the
following fact.

Suppose that a manifold (or finite CW-complex) $M$ is embedded generically in
the space $\R^N$ of a very large dimension, and denote by $M^{*r}$ the union of
all $(r-1)$-dimensional simplices in $\R^N$, whose vertices lie on this
embedded manifold (and the ``genericity'' of the embedding means that if two
such simplices intersect at a certain point, then their minimal faces,
containing this point, coincide). If $M$ is a semialgebraic variety, then by
the Tarski---Seidenberg lemma also $M^{*r}$ is, in particular it has a
triangulation.
\medskip

{\sc Proposition 1} (C.~Caratheodory, see also \cite{Vfil}). {\it The space
$(S^1)^{*r}$ is $PL$-homeomorphic to} $S^{2r-1}$. \quad $\Box$

\section{Proof of the main theorem}

Following \cite{A70}, we use the Alexander duality
\begin{equation}
\label{alex} H^l(HP_d \sm \Sigma_k) \simeq \bar H_{d-l}(\Sigma_k),
\end{equation}
where $\bar H_*$ is the notation for the Borel---Moore homology, i.e. the
homology of the one-point compactification modulo the added point.

To calculate the right---hand group in (\ref{alex}) we construct the {\it
conical resolution} of the space $\Sigma_k$. Let us embed the projective line
$\RP^1$ generically in the space $\R^N$ of a very large dimension, and for any
function $f \in \Sigma_k,$ not equal identically to zero, consider the simplex
$\Delta(f)$ in $\R^N,$ spanned by the images of all points $x_i \in \RP^1,$
corresponding to all possible lines, on which $f$ takes zero value with
multiplicity $\ge k.$ (The maximal possible number of such lines is obviously
equal to $[d/k].$) In the direct product $HP_d \times \R^N$ consider the union
of all simplices of the form $f \times \Delta(f),$ $f \in \Sigma_l \sm 0$. This
union is not closed: the set of its limit points, not belonging to it, is the
product of the point $0 \in HP_d$ and the union of all simplices in $\R^N,$
spanned by the images of no more than $[d/k]$ different points of the line
$\RP^1.$ By the Caratheodory's theorem, the latter union is homeomorphic to the
sphere $S^{2[d/k]-1}.$ We can assume that our embedding $\RP^1 \to \R^N$ is
algebraic, and hence this sphere is semialgebraic. Take a $2[d/k]$-dimensional
semialgebraic disc in $\R^N$ with boundary at this sphere (e.g., the union of
segments connecting the points of this sphere with a generic point in $\R^N$)
and add to the previous union of simplices in $HP_d \times \R^N$ the product of
the point $0 \in HP_d$ and this disc. The obtained set will be denoted by
$\sigma.$ \medskip

{\sc Lemma 2.} {\it The obvious projection $\sigma \to \Sigma_k$ is proper, and
the corresponding map of one-point compactifications of these spaces is a
homotopy equivalence.}
\medskip

This follows easy from the fact that this projection is a stratified map of
semialgebraic spaces, and the preimage of any point $\bar \Sigma_k$ is
contractible, cf. \cite{V}, \cite{Vfil}. \quad $\Box$ \medskip

The space $\sigma$ has a natural increasing filtration: its term $F_p,$ $p \le
[d/k],$ is the union of all $\le (p-1)$-dimensional faces of all simplices,
participating in our construction (or, which is the same, the closure of the
union of all simplices of the form  $f \times \Delta(f)$ over all polynomials
$f$ having no more than $p$ forbidden multiple lines), and $F_{[d/k]+1} =
\sigma$.
\medskip

{\sc Lemma 3.} {\it For any $p =1, \ldots, [d/k],$ the term $F_p \sm F_{p-1}$
of our filtration is the space of a locally trivial fiber bundle over the
configuration space $B(\RP^1,p),$ whose fiber is the direct product of an
$(p-1)$-dimensional open simplex and an $(d+1-pk)$-dimensional real space. The
corresponding bundle of open simplices is orientable if and only if $p$ is odd
(i.e. exactly when the base configuration space is orientable), and the bundle
of $(d+1-pk)$-dimensional spaces is orientable if and only if the number
$k(d+1-pk)$ is even.

The last term $F_{[d/k]+1} \sm F_{[d/k]}$ of this filtration is homeomorphic to
the open $2[d/k]$-dimensional disc.}
\medskip

Indeed, to any configuration $(x_1, \ldots, x_p) \subset \RP^1$ there
corresponds the direct product of the interior part of the simplex in $\R^N,$
spanned by the images of points of this configuration under our embedding
$\RP^1 \to \R^N$, and the subspace in $HP_d,$ consisting of polynomials, having
$k$-fold zeros on corresponding $p$ lines in $\R^2.$ The assertion concerning
the orientations can be checked elementary. \quad $\Box$ \medskip

Consider the spectral sequence $E_{p,q}^r,$ calculating the group $\bar
H_*(\Sigma_k)$ and generated by this filtration. Its term $E_{p,q}^1$ is
canonically isomorphic to the group $\bar H_{p+q}(F_p \sm F_{p-1}).$ \medskip

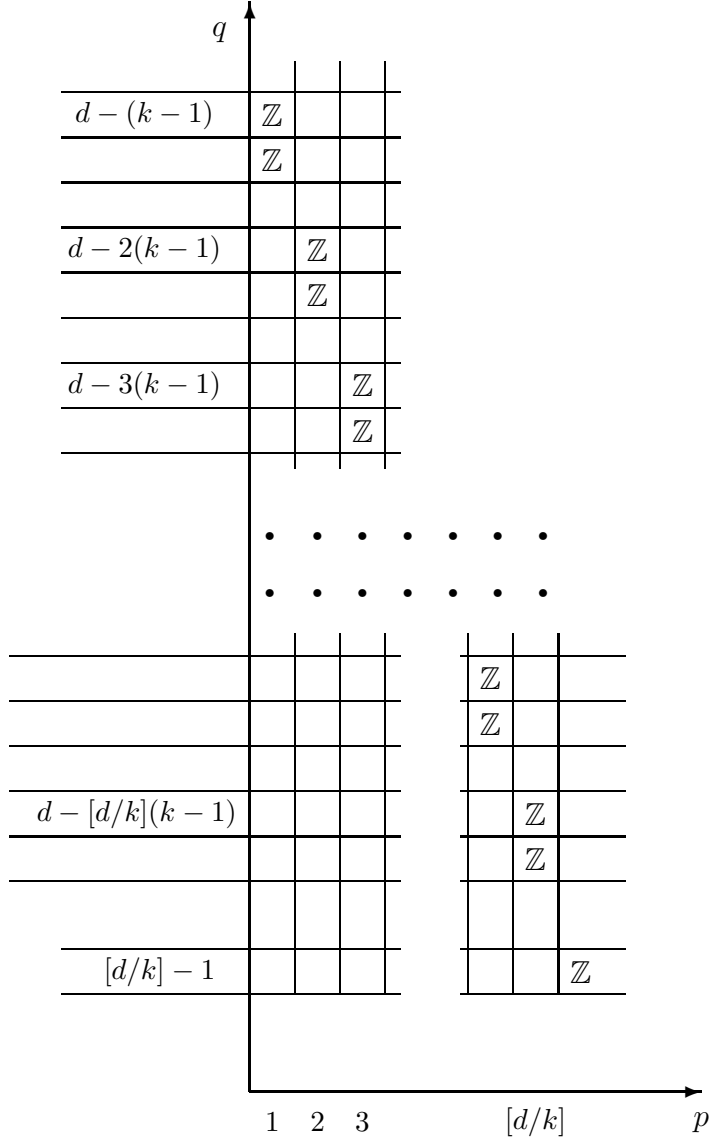
\begin{figure}
\unitlength=1.00mm \special{em:linewidth 0.4pt} \linethickness{0.4pt}
\begin{picture}(100.00,150.00)
\thicklines \put(40.00,10.00){\vector(1,0){60.00}}
\put(40.00,10.00){\vector(0,1){145.00}} \thinlines
\put(40.00,56.00){\line(1,0){20.00}} \put(60.00,62.00){\line(-1,0){20.00}}
\put(40.00,62.00){\line(-1,0){32.00}} \put(8.00,68.00){\line(1,0){52.00}}
\put(43.00,6.00){\makebox(0,0)[cc]{\small $1$}}
\put(55.00,6.00){\makebox(0,0)[cc]{\small $3$}}
\put(49.00,6.00){\makebox(0,0)[cc]{\small $2$}}
\put(69.00,41.00){\line(0,1){30.00}} \put(75.00,71.00){\line(0,-1){30.00}}
\put(81.00,41.00){\line(0,1){30.00}} \put(68.00,68.00){\line(1,0){22.00}}
\put(90.00,62.00){\line(-1,0){22.00}} \put(68.00,56.00){\line(1,0){22.00}}
\put(100.00,6.00){\makebox(0,0)[cc]{$p$}} \put(46.00,111.00){\line(0,1){36.00}}
\put(52.00,147.00){\line(0,-1){36.00}} \put(58.00,111.00){\line(0,1){36.00}}
\put(60.00,143.00){\line(-1,0){45.00}} \put(15.00,137.00){\line(1,0){45.00}}
\put(60.00,131.00){\line(-1,0){45.00}} \put(15.00,125.00){\line(1,0){45.00}}
\put(60.00,119.00){\line(-1,0){45.00}} \put(15.00,113.00){\line(1,0){45.00}}
\put(49.00,122.00){\makebox(0,0)[cc]{${\mathbb Z}$}}
\put(43.00,140.00){\makebox(0,0)[cc]{${\mathbb Z}$}}
\put(36.00,151.00){\makebox(0,0)[cc]{$q$}}
\put(72.00,65.00){\makebox(0,0)[cc]{${\mathbb Z}$}}
\put(42.67,76.33){\circle*{0.94}} \put(49.00,76.33){\circle*{0.94}}
\put(55.00,76.33){\circle*{0.94}} \put(61.00,76.33){\circle*{0.94}}
\put(67.00,76.33){\circle*{0.94}} \put(73.00,76.33){\circle*{0.94}}
\put(79.00,76.33){\circle*{0.94}} \put(79.00,84.00){\circle*{0.94}}
\put(73.00,84.00){\circle*{0.94}} \put(67.00,84.00){\circle*{0.94}}
\put(61.00,84.00){\circle*{0.94}} \put(55.00,84.00){\circle*{0.94}}
\put(49.00,84.00){\circle*{0.94}} \put(42.67,84.00){\circle*{0.94}}
\put(43.00,134.00){\makebox(0,0)[cc]{${\mathbb Z}$}}
\put(46.00,111.00){\line(0,-1){18.00}} \put(52.00,93.00){\line(0,1){18.00}}
\put(58.00,111.00){\line(0,-1){18.00}} \put(60.00,107.00){\line(-1,0){45.00}}
\put(15.00,101.00){\line(1,0){45.00}}
\put(55.00,104.00){\makebox(0,0)[cc]{${\mathbb Z}$}}
\put(15.00,95.00){\line(1,0){45.00}}
\put(55.00,98.00){\makebox(0,0)[cc]{${\mathbb Z}$}}
\put(49.00,116.00){\makebox(0,0)[cc]{${\mathbb Z}$}}
\put(72.00,59.00){\makebox(0,0)[cc]{${\mathbb Z}$}}
\put(60.00,44.00){\line(-1,0){20.00}} \put(40.00,44.00){\line(-1,0){32.00}}
\put(8.00,50.00){\line(1,0){52.00}} \put(68.00,50.00){\line(1,0){22.00}}
\put(90.00,44.00){\line(-1,0){22.00}}
\put(78.00,47.00){\makebox(0,0)[cc]{${\mathbb Z}$}}
\put(40.00,38.00){\line(1,0){20.00}} \put(68.00,38.00){\line(1,0){22.00}}
\put(78.00,41.00){\makebox(0,0)[cc]{${\mathbb Z}$}}
\put(81.00,41.00){\line(0,-1){18.00}} \put(75.00,23.00){\line(0,1){18.00}}
\put(69.00,41.00){\line(0,-1){18.00}} \put(90.00,23.00){\line(-1,0){22.00}}
\put(60.00,23.00){\line(-1,0){45.00}} \put(15.00,38.00){\line(1,0){25.00}}
\put(84.00,26.00){\makebox(0,0)[cc]{${\mathbb Z}$}}
\put(28.00,26.00){\makebox(0,0)[cc]{\small $[d/k]-1$}}
\put(26.00,140.00){\makebox(0,0)[cc]{\small $d-(k-1)$}}
\put(26.00,122.00){\makebox(0,0)[cc]{\small $d-2(k-1)$}}
\put(26.00,104.00){\makebox(0,0)[cc]{\small $d-3(k-1)$}}
\put(25.00,47.00){\makebox(0,0)[cc]{\small $d-[d/k](k-1)$}}
\put(90.00,29.00){\line(-1,0){22.00}} \put(15.00,29.00){\line(1,0){45.00}}
\put(78.00,6.00){\makebox(0,0)[cc]{\small $[d/k]$}}
\put(8.00,38.00){\line(1,0){7.00}} \put(8.00,56.00){\line(1,0){32.00}}
\put(58.00,71.00){\line(0,-1){48.00}} \put(52.00,23.00){\line(0,1){48.00}}
\put(46.00,71.00){\line(0,-1){48.00}}
\end{picture}
\caption{First term of the spectral sequence for even $k$}
\end{figure} \begin{figure}
\unitlength=1.00mm \special{em:linewidth 0.4pt} \linethickness{0.4pt}
\begin{picture}(100.00,135.00)
\thicklines \put(40.00,10.00){\vector(1,0){60.00}}
\put(40.00,10.00){\vector(0,1){130.00}} \thinlines
\put(74.00,41.00){\line(0,1){30.00}} \put(80.00,71.00){\line(0,-1){30.00}}
\put(86.00,41.00){\line(0,1){30.00}} \put(100.00,6.00){\makebox(0,0)[cc]{$p$}}
\put(51.00,111.00){\line(0,1){17.00}} \put(57.00,128.00){\line(0,-1){17.00}}
\put(63.00,111.00){\line(0,1){17.00}}
\put(36.00,137.00){\makebox(0,0)[cc]{$q$}}
\put(77.00,65.00){\makebox(0,0)[cc]{${\mathbb Z}$}}
\put(48.00,76.33){\circle*{0.94}} \put(54.00,76.33){\circle*{0.94}}
\put(60.00,76.33){\circle*{0.94}} \put(66.00,76.33){\circle*{0.94}}
\put(72.00,76.33){\circle*{0.94}} \put(78.00,76.33){\circle*{0.94}}
\put(84.00,76.33){\circle*{0.94}} \put(84.00,84.00){\circle*{0.94}}
\put(78.00,84.00){\circle*{0.94}} \put(72.00,84.00){\circle*{0.94}}
\put(66.00,84.00){\circle*{0.94}} \put(60.00,84.00){\circle*{0.94}}
\put(54.00,84.00){\circle*{0.94}} \put(48.00,84.00){\circle*{0.94}}
\put(51.00,111.00){\line(0,-1){18.00}} \put(57.00,93.00){\line(0,1){18.00}}
\put(63.00,111.00){\line(0,-1){18.00}}
\put(60.00,104.00){\makebox(0,0)[cc]{${\mathbb Z}$}}
\put(60.00,98.00){\makebox(0,0)[cc]{${\mathbb Z}$}}
\put(77.00,59.00){\makebox(0,0)[cc]{${\mathbb Z}$}}
\put(40.00,44.00){\line(-1,0){32.00}} \put(73.00,50.00){\line(1,0){22.00}}
\put(95.00,44.00){\line(-1,0){22.00}} \put(73.00,38.00){\line(1,0){22.00}}
\put(86.00,41.00){\line(0,-1){18.00}} \put(80.00,23.00){\line(0,1){18.00}}
\put(74.00,41.00){\line(0,-1){18.00}} \put(95.00,23.00){\line(-1,0){22.00}}
\put(15.00,38.00){\line(1,0){25.00}}
\put(89.00,26.00){\makebox(0,0)[cc]{${\mathbb Z}$}}
\put(28.00,26.00){\makebox(0,0)[cc]{$[d/k]-1$}}
\put(25.00,47.00){\makebox(0,0)[cc]{\small $d-[d/k](k-1)$}}
\put(95.00,29.00){\line(-1,0){22.00}} \put(83.00,6.00){\makebox(0,0)[cc]{\small
$[d/k]$}} \put(83.00,41.00){\makebox(0,0)[cc]{${{\mathbb Z}_2}$}}
\put(54.00,116.00){\makebox(0,0)[cc]{${{\mathbb Z}_2}$}}
\put(22.00,122.00){\makebox(0,0)[cc]{\small $d-(d-2l)(k-1)$}}
\put(20.00,104.00){\makebox(0,0)[cc]{\small $d-(d-2l+1)(k-1)$}}
\put(8.00,38.00){\line(1,0){7.00}} \put(54.00,7.00){\makebox(0,0)[cc]{\small
$d-2l$}} \put(56.00,0.50){\framebox(17.00,4.50)[cc]{\small $d-2l+1$}}
\put(65.00,5.00){\line(-1,1){5.00}} \put(8.00,56.00){\line(1,0){32.00}}
\put(54.00,122.00){\makebox(0,0)[cc]{$0$}}
\put(83.00,47.00){\makebox(0,0)[cc]{$0$}} \put(14.00,23.00){\line(1,0){26.00}}
\put(50.00,23.00){\line(1,0){15.00}} \put(65.00,29.00){\line(-1,0){15.00}}
\put(40.00,29.00){\line(-1,0){26.00}} \put(50.00,38.00){\line(1,0){15.00}}
\put(65.00,44.00){\line(-1,0){15.00}} \put(50.00,50.00){\line(1,0){15.00}}
\put(65.00,56.00){\line(-1,0){15.00}} \put(73.00,56.00){\line(1,0){22.00}}
\put(95.00,62.00){\line(-1,0){22.00}} \put(63.00,71.00){\line(0,-1){48.00}}
\put(57.00,23.00){\line(0,1){48.00}} \put(51.00,71.00){\line(0,-1){48.00}}
\put(40.00,50.00){\line(-1,0){32.00}} \put(8.00,62.00){\line(1,0){32.00}}
\put(50.00,62.00){\line(1,0){15.00}} \put(73.00,68.00){\line(1,0){22.00}}
\put(65.00,68.00){\line(-1,0){15.00}} \put(40.00,68.00){\line(-1,0){32.00}}
\put(8.00,113.00){\line(1,0){32.00}} \put(50.00,113.00){\line(1,0){15.00}}
\put(65.00,119.00){\line(-1,0){15.00}} \put(40.00,119.00){\line(-1,0){32.00}}
\put(8.00,125.00){\line(1,0){32.00}} \put(50.00,125.00){\line(1,0){15.00}}
\put(65.00,95.00){\line(-1,0){15.00}} \put(40.00,95.00){\line(-1,0){32.00}}
\put(8.00,101.00){\line(1,0){32.00}} \put(50.00,101.00){\line(1,0){15.00}}
\put(65.00,107.00){\line(-1,0){15.00}} \put(40.00,107.00){\line(-1,0){32.00}}
\end{picture}
\caption{Term $E^1$ for odd $k$}
\end{figure}
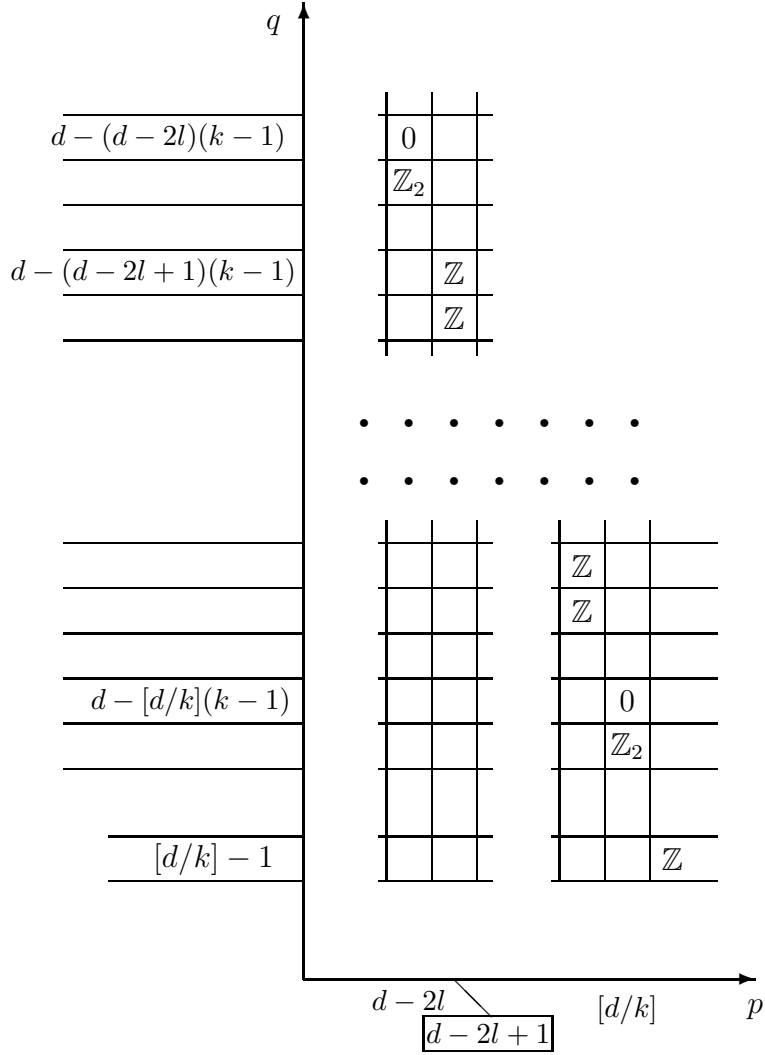

{\sc Corollary.} {\it The term $E_{p,q}^1$ of our spectral sequence is as in
Fig. 1 for even $k$ and as in Fig. 2 is $k$ is odd, i.e.:

if $k$ is even or $d-p$ is odd, $1 \le p \le [d/k],$ then in the column
$E^1_{p,*}$ only the cells $E^1_{p,q}$ with $q=d-p(k-1)$ and $q=d-p(k-1)-1$ are
nontrivial; these two cells are isomorphic to $\Z$;

if the number $k(d-p-1)$ is odd, $1 \le p \le [d/k],$ then the column
$E^1_{p,*}$ contains unique nontrivial cell $E^1_{p,d-p(k-1)-1} \simeq \Z_2$;

for any $k$ and $d$ the unique nontrivial cell of the column $E^1_{[d/k]+1,*}$
is $E^1_{[d/k]+1,[d/k]-1} \simeq \Z;$

all other columns of the spectral sequence are trivial.} \quad $\Box$ \medskip

{\sc Remark.} The situation shown in columns $p=[d/k]$ and $p=[d/k]-1$ of Fig.
2 appears when $d-[d/k]$ is even. In fact the opposite also can happen, but in
the most interesting case, when $d$ is a multiple of $k,$ the situation is
exactly as in the picture. In this last case the unique nontrivial cell
$E^1_{[d/k]+1,[d/k]-1} \simeq \Z$ of the column $p=[d/k]+1$ lies in the same
horizontal row as the unique nontrivial cell $E^1_{[d/k],[d/k]-1} \simeq \Z_2$
of the column $p=[d/k].$
\medskip

{\sc Proposition 2.} {\it For any $d, k,$ except for the case when $k$ is odd
and $d$ is a multiple of $k$, our spectral sequence degenerates in the term
$E_1$, i.e. $E_1 = E_\infty$. In the exceptional case there is unique
nontrivial operator $d_1,$ acting from the cell $E^1_{d/k+1,d/k-1} \sim \Z$ to
the cell $E^1_{d/k,d/k-1} \sim \Z_2,$ killing the latter cell, after which the
spectral sequence also degenerates.} \medskip

Main theorem follows immediately from this proposition.
\medskip

{\it Proof of Proposition 2.} For almost all $d$ and $k$ the assertion follows
immediately from the explicit form of the term $E_1.$ The only two cases, when
nontrivial differentials could occur, are the case $k=2$ (a half of which is
already studied, see Main example in \S \ 2, and the remaining case of odd $d$
is even so easy) and the case when $d$ is a multiple of $k$. In the latter case
(for $k>2$) the unique nontrivial differential can act from the cell
$E^1_{d/k+1,d/k-1} \sim \Z$ to the cell $E^1_{d/k,d/k-1}.$ If $k$ is even, then
the latter cell is isomorphic to $\Z,$ and the corresponding term $F_{d/k} \sm
F_{d/k-1}$ of our filtration is the oriented bundle with fiber $\R^1,$ whose
base also is oriented and is, in its turn, the space of a fibre bundle with
base $B(S^1,d/k)$ and the fiber equal to an $(d/k-1)$-dimensional open simplex.
The boundary of the disc $F_{d/k+1} \sm F_{d/k}$ in this term of filtration
coincides with the zero section of the former (line) bundle, therefore it is
the boundary (modulo lower terms of the filtration) of the space of the bundle
with fiber $\R_+$.

In the case of odd $k$ similar line bundle is non-orientable, therefore when we
try to span this zero section by a similar chain over a maximal
simple-connected domain in the base, then we construct a homology between the
image of this section and a cycle generating the group $H_{2d/k-1}(F_{d/k} \sm
F_{d/k-1}) \sim \Z_2$. (An adequate picture here is the M\"obius band, whose
equator circle is the boundary of the disc $F_{d/k+1} \sm F_{d/k}$.) All the
further differentials are trivial by the dimensional reasons, and main theorem
is completely proved.
\medskip

{\sc Acknowledgment.} I thank A.~Durfee for a stimulating discussion and
G.~Kalai, who indicated me that Proposition 1 is due to Caratheodory.


\begin{thebibliography}{99}

\bibitem{A70} V.~I.~Arnold, {\it On some topological invariants of
algebraic functions,} Trans. Moscow Math. Soc. 1970, V.~21, p.~30--52.

\bibitem{A89} V.~I.~Arnold, {\it Spaces of functions with
moderate singularities,} Funct. Anal. and its Appl. 23:3, (1989), p.~1--10.

\bibitem{V89} V.~A.~Vassiliev, {\it Topology of spaces of functions without
complicated singularities.} Funct. Anal. and its Appl., 1989, 23:4, p.~24--36.

\bibitem{V90} V.~A.~Vassiliev, {\it Cohomology of knot spaces,}
In: Adv. in Sov. Math.; Theory of Singularities and its Appl. (ed.
V.~I.~Arnold). AMS, Providence, R.I., 1990, p.~23--69.

\bibitem{V91} V.~A.~Vassiliev, {\it Geometric realization of the homology of
classical Lie groups and complexes, $S$-dual to flag manifolds}, St.-Petersburg
Math. J. 3:4 (1991), p.~108--115.

\bibitem{V} V.~A.~Vassiliev, {\it Complements of discriminants of smooth
maps: topology and applications. Revised ed.}, Transl. of Math. Monographs,
v.~98, AMS, Providence RI., 1994.

\bibitem{Vfil} V.~A.~Vassiliev, {\it Topology of complements of discriminants,}
Moscow,  Phasis, 1997, 552 p. (in Russian).

\end{thebibliography}
\end{document}